\documentclass[twocolumn,secnumarabic,amssymb, aps, prd]{revtex4-2}  

\usepackage{amssymb,amsmath,amsthm}
\usepackage{mathtools}
\usepackage{amscd}
\usepackage{enumerate}
\usepackage{subcaption}
\usepackage{caption}
\usepackage{graphicx}

\setcitestyle{authoryear,round}
\usepackage{color}

\usepackage{newtxtext}

\numberwithin{equation}{section} \oddsidemargin=-.0cm
\usepackage[font=small,labelfont=bf]{caption}

\usepackage{color}
\usepackage{graphicx,graphics}
\usepackage{amsfonts}

\usepackage{sidecap}
\usepackage{dsfont}

\usepackage{bm} 

\sidecaptionvpos{figure}{c}

\usepackage{xcolor}
\usepackage[colorlinks=true, pdfstartview=FitV, linkcolor=blue,
            citecolor=blue, urlcolor=blue]{hyperref}  

\usepackage{sidecap}
\graphicspath{{./figures/},{./Figures/}}

\newtheorem{thm}{Theorem}[section]
\newtheorem{lem}{Lemma}[section]
\newtheorem{defi}{Definition}[section]

\newtheorem{prop}[thm]{Proposition}

\newtheorem{rem}{Remark}[section]
\newtheorem{cor}{Corollary}[section]

\def\bt{\begin{thm}}
\def\et{\end{thm}}
\def\bl{\begin{lem}}
\def\el{\end{lem}}
\def\bd{\begin{defi}}
\def\ed{\end{defi}}
\def\bc{\begin{cor}}
\def\ec{\end{cor}}
\def\bp{\begin{proof}}
\def\ep{\end{proof}}
\def\br{\begin{rem}}
\def\er{\end{rem}}

\def\bprop{\begin{prop}}
\def\eprop{\end{prop}}

\def\d{\, \mathrm{d}}

\def\be{\begin{equation}}
\def\ee{\end{equation}}
\def\bes{\begin{equation*}}
\def\ees{\end{equation*}}
\def\bea{\begin{equation} \begin{aligned}}
\def\eea{\end{aligned} \end{equation}}
\def\beas{\begin{equation*} \begin{aligned}}
\def\eeas{\end{aligned} \end{equation*}}

\def\bpp{\begin{proof}}
\def\epp{\end{proof}}

\def\bi{\begin{itemize}}
\def\ei{\end{itemize}}
\def\ben{\begin{enumerate}}
\def\een{\end{enumerate}}

\def\y{\boldsymbol{y}}

\definecolor{rred}{rgb}{0.7,0,0.1}

\definecolor{greenrb}{rgb}{0.2,0.6,0.2}

\def\cH{\mathcal H}

\begin{document}

\title{Galerkin-Koornwinder approximations of   delay differential equations for physicists}

\author{Micka\"el D. Chekroun}
\email{mchekroun@atmos.ucla.edu} 
\affiliation{Department of Earth and Planetary Sciences, Weizmann Institute of Science, Rehovot 76100, Israel,}
\affiliation{Department of Atmospheric and Oceanic Sciences, University of California, Los Angeles,  CA 90095-1565, USA}

\author{Honghu Liu}
\email{hhliu@vt.edu}
\affiliation{Department of Mathematics, Virginia Tech, Blacksburg, VA 24061, USA}

\begin{abstract}
Formulas for Galerkin-Koornwinder (GK) approximations of delay differential equations are summarized.
The functional analysis ingredients (semigroups, operators, etc.) are intentionally omitted to focus instead on the formulas required to perform GK approximations in practice. 

\end{abstract}

\maketitle


\section{Galerkin-Koornwinder (GK) approximations of DDEs}

We consider nonlinear {\it scalar} differential delay equations (DDEs) of the form
\bea \label{Eq_DDE}
\frac{\d x(t)}{\d t}& = a x(t) + b x(t-\tau) + c \int_{t-\tau}^t x(s)\d s \\
 &\;\; + F\Big(x(t), x(t-\tau),\int_{t-\tau}^t x(s) \d s\Big),
\eea
where $a$, $b$ and $c$ are real numbers, $\tau> 0$ is the delay parameter, and $F$
is a nonlinear function. 
We restrict ourselves to the scalar case to simplify the maths, but the approach extends 
to systems of nonlinear DDEs involving possibly several delays as detailed in \citep{CGLW16} and illustrated in \citep{CKL20} for the cloud-rain delay model of \citep{koren2011aerosol}.

\subsection{Koornwinder polynomials}
First, let us recall that Koornwinder polynomials $K_n$ are obtained from 
Legendre polynomials $L_n$, for any nonnegative integer $n$, according to the relation
\bea \label{eq:Pn}
&K_n(s)= -(1+s)\frac{\d}{\d s} L_n(s) +( n^2 + n + 1) L_n(s),\\
& \; s \in [-1, 1];
\eea
see \cite[Eq.~(3.3)]{CGLW16}.

Koornwinder polynomials are  known to form an orthogonal set for the following weighted inner product on $[-1,1]$ with a point-mass, $\mu(\d x)= \frac{1}{2} \d x  + \delta_{1},$
where $\delta_1$ denotes the Dirac point-mass at the right endpoint $x=1$; see  \cite{Koo84}.  In other words, the following orthogonality property holds:
\begin{widetext}
\bea
\int_{-1}^{1} K_n(s) K_p(s) \d \mu (s)& = \frac{1}{2} \int_{-1}^{1} K_n(s) K_p(s) \d s + K_n(1) K_p(1)\\
&=0, \, \mbox{ if } p\neq n.
\eea
\end{widetext}

It is also worthwhile noting that Koornwinder polynomials augmented by the right endpoint values as follows 
\be \label{eq:Pn_prod}
\mathcal{K}_n = (K_n, K_n(1)),
\ee
form an orthogonal basis of the  product space 
\bes \label{eq:E}
\mathcal{E} = L^2([-1,1); \mathbb{R}) \times  \mathbb{R},
\ees
endowed with the inner product:
\be \label{eq:inner_E}
\Big\langle (f, a), (g, b) \Big\rangle_{\mathcal{E}}  = \frac{1}{2} \int_{-1}^1 f(s)g(s) \d s  + ab.
\ee
The norm induced by this inner product is denoted by $\|\cdot\|_{\mathcal{E}}$. The basis function $\mathcal{K}_n$ has then its $\|\cdot\|_{\mathcal{E}}$-norm  given by \cite[Prop.~3.1]{CGLW16}:
\be \label{eq:Kn_norm}
\|\mathcal{K}_n\|_{\mathcal{E}} = \sqrt{\frac{(n^2+1)((n+1)^2+1)}{2n+1}}.
\ee
This is a useful property for calculating the GK approximations, as it will come apparent below.

Finally, since the original Koornwinder basis is given on the interval $[-1, 1]$ and the state space of a DDE such as Eq.~\eqref{Eq_DDE} is made of functions defined on $(-\tau,0)$, we will work mainly with the following rescaled version of  Koornwinder polynomials. 
The rescaled Koornwinder polynomials $K_n^\tau$ are defined as follows 
\bea \label{eq:Pn_tilde}
K^\tau_n\colon  [-\tau, 0] & \rightarrow \mathbb{R}, \\
\theta & \mapsto  K_n \Bigl( 1 + \frac{2 \theta }{\tau} \Bigr).
\eea
They form orthogonal polynomials on the interval $[-\tau, 0]$ for the $L^2$-inner product on $(-\tau,0)$ with a Dirac point-mass at the right endpoint $0$.

The following family of rescaled Koornwinder polynomials augmented with a right endpoint value,
\be \label{eq:Pn_tilde_prod}
\mathcal{K}_n^\tau = (K_n^\tau, K_n^\tau(0)),
\ee
forms then an orthogonal basis for the Hilbert space $\mathcal{H} = L^2([-\tau,0); \mathbb{R}) \times  \mathbb{R}$ endowed with the inner product $\langle \cdot, \cdot \rangle_{\mathcal{H}}$ given by \eqref{eq:inner_E} in which the integral is taken on $(-\tau,0)$ and weighted by $1/\tau$. Note that since $K_n(1)=1$ \cite[Prop.~3.1]{CGLW16}, we have also 
\be\label{Eq_normalization}
K_n^\tau(0)  = 1.
\ee
Finally, observe that $\|\mathcal{K}_n^\tau\|_{\cH} = \|\mathcal{K}_n\|_{\mathcal{E}}$.

\subsection{Space-time representation and GK approximations}
Unlike ordinary differential equations (ODEs), a DDE such as Eq.~\eqref{Eq_DDE} is an infinite-dimensional dynamical system in which a history segment has to be specified over the interval $[-\tau,0)$ to apprehend properly the existence and computation problems of its solutions. 
  
\begin{figure}[hbtp]
   \centering
\includegraphics[width=0.5\textwidth, height=0.5\textwidth]{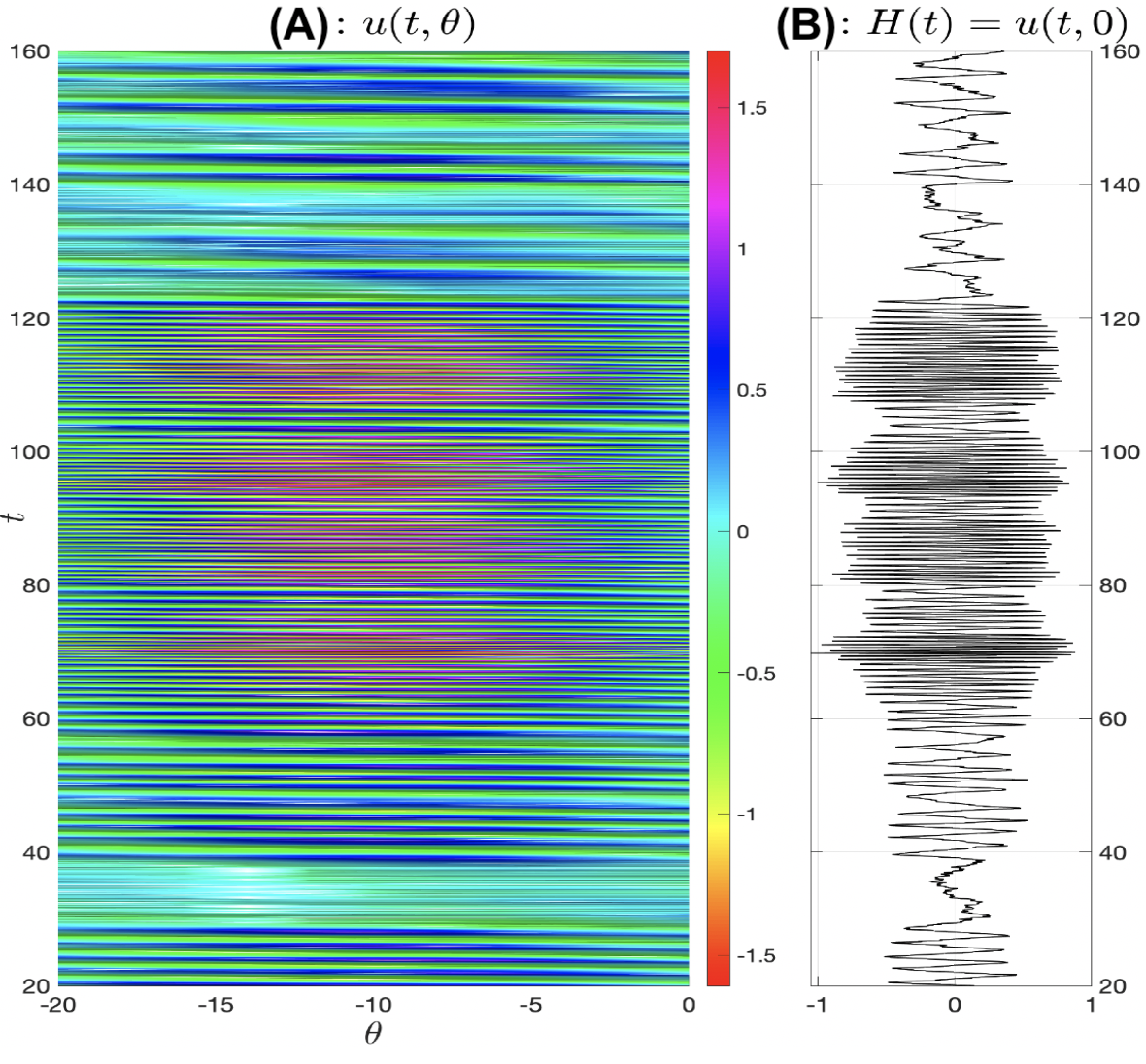}
  \caption{{\small  Panel A shows a solution $u(t, \theta)$ to a transport problem associated with a DDE. In this example, the boundary value $u(t, 0)$ at $\theta=0$ is shown in panel B as the time series $H(t)$ (from the Supplementary Material of \citep{Chekroun_al22SciAdv}). This time series at the right endpoint $\theta=0$, coincides, in the deterministic setting, with the time series $x(t)$ obtained by solving the original DDE.}}   \label{Fig_hovmoller}
\end{figure}

Thus, given a function of time, $x$, solving Eq.~\eqref{Eq_DDE} one distinguishes between the {\it history segment} $\{x(t + \theta): \theta \in [-\tau, 0)\}$ and the {\it current state}, $x(t)$. Denoting by $u(t,\theta)$ the history segment, one can then rewrite  the DDE  \eqref{Eq_DDE} as the transport equation 
\be\label{lin_PDE}
\partial_t u = \partial_{\theta}u, \quad -\tau \le \theta < 0,
\ee
subject to the {\it nonlocal and nonlinear boundary condition}
\bea\label{PDE_BC}
\partial_{\theta} u|_{\theta = 0}& = a u(t,0) + b u(t,-\tau) + c \int_{t-\tau}^t u(s,0)\d s\\
&\;\; + F\Big(u(t,0), u(t,-\tau),\int_{t-\tau}^t  u(s,0) \d s\Big),
\eea
for $t \geq 0$.
This reformulation is often called the space-time representation in the literature; see Fig.~\ref{Fig_hovmoller} for an illustration.

Using rescaled  Koornwinder polynomials \eqref{eq:Pn_tilde}, one can then seek for  approximations of $u(t,\theta)$ solving  Eqns.~\eqref{lin_PDE}-\eqref{PDE_BC} under the form
\be\label{Eq_App1}
u_N(t,\theta) = \sum_{j=0}^{N-1} y_j(t) K_j^\tau (\theta).
\ee
Given the property \eqref{Eq_normalization}, the current state $x(t)=u(t,0)$ is approximated 
by
\be\label{Eq_App2}
x_N(t)= \sum_{j=0}^{N-1} y_j(t). 
\ee
Note that in \eqref{Eq_App1} the Koornwinder polynomials are indexed according to their degree $j$.  

The question arises then of determining the equations that the coefficients $y_j(t)$ in \eqref{Eq_App1} must satisfy in order to have that  $u_N(t,\theta)$ provides a rigorous approximation that converges to $u(t,\theta)$ as $N\rightarrow \infty.$
This nontrivial problem has been solved in \citep{CGLW16}. The next section recalls the structure of these equations forming what we call a GK system.

\vspace{.5cm}
\subsection{The GK approach in action}
The analysis of  \citep{CGLW16} shows that $y_j(t)$ solving the $N$-dimensional system of ODEs \eqref{Galerkin_AnalForm} provides a rigorous method for approximating the solutions to a broad class of DDEs. 
Such systems are given by: 
\begin{widetext}
\begin{equation} \label{Galerkin_AnalForm}
\begin{aligned}
\frac{\d y_j}{\d t} & = \frac{1}{\|\mathcal{K}_j\|_{\mathcal{E}}^2 } \sum_{n=0}^{N-1} \bigg( a + b K_n(-1) + c \tau (2 \delta_{n,0} - 1) \\
& \hspace{8em} + \frac{2}{\tau}\sum_{k=0}^{n-1} a_{n,k} \left( \delta_{j,k} \|\mathcal{K}_j\|^2_{\mathcal{E}} - 1 \right) \bigg) y_n(t) \\
& \hspace{1em} + \frac{1}{\|\mathcal{K}_j\|_{\mathcal{E}}^2} F \left( \sum_{n=0}^{N-1} y_n(t),  \sum_{n=0}^{N-1} y_n(t) K_n(-1),\tau y_0(t) - \tau \sum_{n=1}^{N-1} y_n(t) \right),\\
& \;\; 0\leq j\leq N-1.
\end{aligned}
\end{equation}
\end{widetext}
Here the Kronecker symbol $\delta_{j,k}$ has been used, and the coefficients $a_{n,k}$ are obtained by solving a triangular linear system in which the right-hand side has 
explicit coefficients depending on $n$ \cite[Prop.~5.1]{CGLW16}; see Appendix~\ref{sect:coef_matrix_proof}. 
The values of the $K_n(-1)$'s are known and recalled in \eqref{Eq_K-1} below.


\begin{widetext}
We rewrite now the GK system \eqref{Galerkin_AnalForm} into the following compact form:
\be \label{Galerkin_cptForm}
\boxed{\frac{\d \y}{\d t} = A_N (\tau) \y + F_N (\y,\tau),}
\ee
where $\y= (y_0, \cdots, y_{N-1})^{T}$, 
and in which the $A_N (\tau) \y$-term and $F_N(\y,\tau)$-term group the linear and nonlinear terms in Eq.~\eqref{Galerkin_AnalForm}, respectively. 
The entries of the $N\times N$ matrix  $A_N (\tau)$ are given by 
\begin{equation} \label{eq:A}
\boxed{
\begin{aligned}
\Big(A_N (\tau)\Big)_{i,j}& = \frac{1}{\|\mathcal{K}_i\|_{\mathcal{E}}^2 }\Big(a + b K_j(-1) + c \tau (2 \delta_{j,0} - 1) \\
& \hspace{8em}+ \frac{2}{\tau}\sum_{k=0}^{j-1} a_{j,k} \left( \delta_{i,k} \|\mathcal{K}_i\|^2_{\mathcal{E}} - 1 \right ) \Big), \;\; i, j = 0, \cdots, N-1.
\end{aligned}
}
\end{equation}
 Here, the coefficients $a_{j,k}$ are independent of the DDE model and are determined by \cite[Prop.~5.1]{CGLW16}; see Proposition \ref{prop:dPn} below. 
Re-writing $A_N(\tau)$ as $A_N(\tau)=2/\tau P_N +Q_N$, one observes that only $Q_N$ depends on the model's  parameters. More generally, the matrix $A_N(\tau)$  accounts for the approximation of the transport equation $\partial_t u=\partial_\theta u$ and the contribution of the linear terms arising in the nonlocal boundary condition \eqref{PDE_BC}; see \citep{CGLW16}.

The nonlinear mapping $F_N$ has its components $F_{N}^j$ given by 
\be \label{eq:G}
\boxed{F_{N}^j(\y,\tau) = \frac{1}{\|\mathcal{K}_{j}\|_{\mathcal{E}}^2} F \left( \sum_{n=0}^{N-1} y_n(t),  \sum_{n=0}^{N-1} y_n(t) K_n(-1),\tau y_0(t) - \tau \sum_{n=1}^{N-1} y_n(t) \right),}
\ee
for each  $0\leq j\leq N-1$.
\end{widetext}

Due to rigorous convergence results of GK approximations \citep{CGLW16}, the GK formulas above provide a powerful apparatus to analyze DDEs by means of ODE approximations.

\subsection{Example: GK  approximations of the Suarez and Schopf model}
The Suarez and Schopf model  is given by the following DDE  \citep{Suarez_al88} 
\be \label{eq:SS}
\frac{\d T}{\d t} = T(t) - \alpha T(t - \tau) - T^3(t),
\ee
where $\tau$ and $\alpha$ are positive constants; and the physically relevant range of $\alpha$ used in \citep{Suarez_al88} is $(0,1)$. We refer to \citep{Suarez_al88} for the physical interpretation of this model in relationship with the El Ni\~no-Southern Oscillation (ENSO) phenomenon, and to \citep{CGN18} and references therein for other ENSO models with delays.  

It is clear that for this given range of $\alpha$, Eq.~\eqref{eq:SS} admits three fixed points: 
\bes
T_0 = 0, \qquad T_{+} = \sqrt{1- \alpha}, \qquad T_{-} = - \sqrt{1- \alpha}.
\ees

Let us first introduce the perturbed variable 
\be
x = T - T_{+}.
\ee 
Equation \eqref{eq:SS}, when written in the perturbed variable $x$, reads
\bea \label{Eq_SS_perturb}
\frac{\d x}{\d t}& = (1 - 3 T_{+}^2) x(t) \\
&\qquad- \alpha x(t - \tau)  - 3 T_{+} x^2(t) - x^3(t).
\eea
The above equation fits into Eq.~\eqref{Eq_DDE} with 
\bea \label{eq:SS_pars}
&a = 1 - 3 T_{+}^2, \quad b = -\alpha,  \quad c = 0, \\
& \text{and} \quad F(x(t)) = - 3 T_{+} x^2(t) - x^3(t).
\eea 
The corresponding $N$-dimensional GK approximation can then be obtained by using Eqns.~\eqref{Galerkin_cptForm}--\eqref{eq:G}, and it reads 
\be  \label{Eq_Galerkin_SS}
\frac{\d \y}{\d t} = A(\tau) \y  +  G (\y),
\ee
dropping the dependence on $N$.

The entries of the matrix $A(\tau)$ are given by
\bea \label{eq:A_SS}
(A(\tau))_{j,n} &  = \frac{1}{\|\mathcal{K}_j\|_{\mathcal{E}}^2 } \sum_{n=0}^{N-1}  \Bigl( 1 - 3 T_{+}^2  - \alpha K_n(-1)  \\
&  \qquad +  \frac{2}{\tau}\sum_{k=0}^{n-1} a_{n,k} \left( \delta_{j,k} \|\mathcal{K}_j\|^2_{\mathcal{E}} - 1 \right ) \Bigr), 
\eea
for $j, n  = 0, \cdots, N-1$. The $ a_{n,k}$ are determined thanks to Proposition \ref{prop:dPn}. 

The nonlinearity $G$ is given as the following sum of $N$-dimensional mappings
\bes \label{eq:SS_G}
 G(\y) =  G_2(\y) +  G_3(\y),
\ees
with 
\beas \label{Eq_G_decomp_SS}
G_2(\y) &= - 3 T_{+}  \left( \sum_{n=0}^{N-1} y_n \right)^2 \boldsymbol{\nu}_N, \\
G_3(\y) &= - \left( \sum_{n=0}^{N-1} y_n \right)^3 \boldsymbol{\nu}_N,
\eeas
 where $\boldsymbol{\nu}_N$ denotes the $N$-dimensional column vector  
 \be\label{Eq_vector_norms}
 \boldsymbol{\nu}_N=\Big(\frac{1}{\|\mathcal{K}_0\|_{\mathcal{E}}^2}, \cdots, \frac{1}{\|\mathcal{K}_{N-1}\|_{\mathcal{E}}^2} \Big)^T. 
\ee
Note that here only the linear terms in the GK-approximations depend on $\tau$ as the DDE model \eqref{Eq_SS_perturb} contains delay terms only in its linear ones. If delay terms would be present in the nonlinear terms as well then $G$ would depend on $\tau$ according to  \eqref{eq:G}.

We set $\alpha = 0.75$. The steady state is then given by $T_{+} = \sqrt{1 - \alpha} = 0.5$.
Setting $N=6$, the 6-D GK approximation system of the form \eqref{Eq_Galerkin_SS} reads here as:
\bea  \label{Eq_Galerkin_SS_6D}
&\frac{\d \y}{\d t} = \left(M_1 + \frac{1}{\tau} M_2 \right) \y \\
 &\;\; - \left(\sum_{n=1}^{6} y_n \right)^2  \left( 3 T_{+}  + \sum_{n=1}^{6} y_n \right) \boldsymbol{\nu}_6,
\eea
where $\y= (y_1, \cdots, y_6)^{T}$, $\boldsymbol{\nu}_6$ is given by \eqref{Eq_vector_norms} with $N=6$, and the matrix $A(\tau)$ in \eqref{Eq_Galerkin_SS} is given by $M_1 + \frac{1}{\tau} M_2$ with 
\begin{widetext}
\be \label{Eq_M1}
M_1 = \begin{pmatrix}
-0.25  &  1.25 &   -2.5 &    5 &   -7.75 &  11.75 \\
   -0.15 &   0.75 &   -1.5 &  3 &  -4.65 &  7.05  \\
   -0.05 &    0.25 &   -0.5 &  1 &   -1.55 &    2.35 \\
   -0.0206  &   0.1029 &   -0.2059 &    0.4118  &  -0.6382 &    0.9676 \\
   -0.0102  &  0.0509  & -0.1018  &  0.2036 &  -0.3156  &  0.4785 \\
   -0.0057  &  0.0286 &  -0.0572  &  0.1143 &  -0.1772  &  0.2687
\end{pmatrix},
\ee
\be \label{Eq_M2}
M_2 = \begin{pmatrix}
         0 &   2   &  -3 &    7 &  -10 &   16  \\
         0 &  -1.2 &   7.8 &  -10.2 &   20.4 &  -26.4 \\
         0 &  -0.4 &    -2.4 &   11.6 &  -12.2 &   24.2 \\
         0 &  -0.1647 &   -0.9882 &   -3.4588 &  14.7412 &  -13.0941 \\
         0 &  -0.0814 &  -0.4887 &  -1.7104 &   -4.4796  & 17.7557 \\
         0 &  -0.0457 &  -0.2744  & -0.9605  &  -2.5156  & -5.4886 
\end{pmatrix}.
\ee
We refer to Appendix~\ref{sect:coef_matrix_proof} for the computation of  $M_1$ and $M_2.$
\end{widetext}

Numerical results show that the GK system \eqref{Eq_Galerkin_SS_6D} with  $N=6$ exhibits already very good skills in approximating the solutions to the DDE model \eqref{Eq_SS_perturb}; see Fig.~2 in Supplementary Material of \citep{CLM20_closure}. In general, a few Koornwinder polynomials are needed to capture the dynamics in absence of time-dependent forcings; see also Fig.~4 in \citep{CGLW16} and \citep[Sec.~4.3.3]{CKL18_DDE}.

\appendix

\section{Determination of $M_1$ and $M_2$ in Eq.~\eqref{Eq_Galerkin_SS_6D}} \label{sect:coef_matrix_proof}
The calculations of $M_1$ and $M_2$ are based on the RHS of \eqref{eq:A_SS}. First, the terms 
$$\frac{1}{\|\mathcal{K}_j\|_{\mathcal{E}}^2 } \sum_{n=0}^{N-1}  ( 1 - 3 T_{+}^2  - \alpha K_n(-1)),$$ when collected for each $j$ and $n$ allow us to build up $M_1$ readily once  one recalls \eqref{eq:Kn_norm} and that due to \eqref{eq:Pn},
\be\label{Eq_K-1}
K_n(-1)=( n^2 + n + 1) (-1)^n,
\ee
since $L_n(-1)=(-1)^n$ from the properties of the Legendre polynomials.

The determination of $M_2$ which results from the collection of the terms $2 \sum_{k=0}^{n-1} a_{n,k} \left( \delta_{j,k} \|\mathcal{K}_j\|^2_{\mathcal{E}} - 1 \right )$, requires  the knowledge of the coefficients $a_{n,k}$. 
These coefficients arise in the expression of the derivatives of the Koornwinder polynomials in terms of these polynomials themselves; see \cite[Appendix~B]{CGLW16} for more details.  

The coefficients $a_{n,k}$ are then obtained by simply solving the triangular system \eqref{eq:algebraic} described in the following proposition.

\begin{prop} \label{prop:dPn}
The Koornwinder polynomial $K_n$ of degree $n$ defined in \eqref{eq:Pn} satisfies the differential relation
\be \label{eq:dPn}
\frac{\d K_n}{\d s}(s) = \sum_{k = 0}^{n-1} a_{n,k} K_k(s), \quad s \in (-1,1),
\ee
where the vector $\boldsymbol{a}_n$ made of the  $a_{n,k}$-coefficients,  
$$\boldsymbol{a}_n=(a_{n,0}, \cdots, a_{n,n-1})^T,$$
 solves the  following upper triangular system
\be \label{eq:algebraic}
\mathbf{T}\boldsymbol{a}_n = \boldsymbol{b}_n,
\ee
with $\mathbf{T}=(\mathbf{T}_{i,j})_{n\times n}$ and $\boldsymbol{b}_{n}=(b_{n,0}, \cdots, b_{n,n-1})^T$ given by
\begin{widetext}
\bea \label{eq:algebraic_def}
\mathbf{T}_{i,j} & = \begin{cases}
0, & \; \text{ if } j < i,\\
i^2 + 1, & \; \text{ if } j = i,\\
-(2i+1), & \; \text{ if } j > i,
\end{cases} \; \quad \text{ where } \quad 0 \le i, j \le n-1, \\
b_{n,i} & = \begin{cases}
 -\frac{1}{2}(2i+1)(n+i+1)(n-i), & \text{ if $n+i$ is even}, \vspace{1em}\\
(n^2 + n)(2i+1) - \frac{i}{2}(n+i)(n-i+1) & \\
\hspace{2em} -\frac{1}{2}(i+1)(n-i-1)(n+i+2), & \text{ if $n+i$ is odd}.
\end{cases}
\eea
\end{widetext}

Finally, the rescaled Koornwinder polynomials  satisfy 
\be \label{eq:dKn_tau}
\frac{\d K^\tau_n}{\d \theta}(\theta) = \frac{2}{\tau} \sum_{k = 0}^{n-1} a_{n,k} K^\tau_k(\theta),  \quad \theta\in (-\tau,0).
\ee

\begin{rem}
We note that the formula for $b_{n,i}$ that appeared in \citep{CGLW16} and in \citep{CKL20} contains two typos (a sign error and an extra factor $i$ in one of the terms) that are here rectified in the formula \eqref{eq:algebraic_def}. We would like to express our gratitude to Tetsushi Saburi for pointing out these typos to us. We emphasize though that it is only a typo as the correct formulas were implemented in our codes, thus not questioning the diverse numerical results published in \citep{CGLW16}, \citep{CKL18_DDE}, \citep{CKL20}, \citep{Chekroun_al22SciAdv}, and in the Supplementary Material of \citep{CLM20_closure}.
\end{rem}

\end{prop}



\bibliography{reference}

\end{document}